\documentstyle[12pt,amstex]{article}
\begin{document}
\newcommand{\R}{\hbox{I\!R}}
\newcommand{\Hi}{\hbox{I\!H}}
\newcommand{\tp}{\mbox{$\tilde p$}}
\newcommand{\wm}{\mbox{$\widetilde M$}}
\newcommand{\wu}{\mbox{$\widetilde U$}}
\begin{center}

\centerline{\Large On the number of geodesic segments connecting two points}

\vspace{0.25in}

{\Large  on manifolds of non-positive curvature}

\vspace{0.5in}

by

\vspace{0.25in}

Paul Horja

\end{center}

\vspace{0.5in}

{\bf 1. Introduction}

Let $M$ be a compact manifold Riemannian manifold of dimension $n \geq 2$,
with a metric of sectional curvature bounded above by $\chi \leq 0$
(non-positive curvature). In this paper we prove that in the case
of negative curvature ($\chi < 0$) on such manifolds
there exist pairs of points connected by at least $2n+1$ geometrically
 distinct geodesic segments (i.e. length minimizing). A class of points which
provide examples in this class are the points situated at distance equal to the
diameter of the manifold. A simplified version of the  method allows us to show
 that in the case of non-positive curvature ($\chi = 0$) for any point there
exist another point and $n+1$ geometrically distinct geodesic segments
connecting them. The essential ingredient in the proofs is the basic metric
 property of the spaces of non-positive and negative curvature to have their
distance function convex and, in a sense that will be explained in the paper,
even almost strictly convex. The results can also be seen as estimates for
 the ``order" of the points in the cut locus for manifolds of non-positive
 curvature.

In the case of positive curvature the situation changes: for the ellipsoid
in $\R ^3$ with axes of different lengths, the points at maximal
distance are connected by two geodesic segments, but for the sphere by
infinitely many geodesic segments. For the flat torus obtained as a
quotient of $\R^2$ by a lattice not generated by two orthogonal vectors,
the maximal ``order'' of the points in the cut locus is $3$. Interesting is
the situation for convex polyhedra in $\R^3$  which is intermediate between
the cases of negative and positive curvature. For a large class of them,
namely for those
admitting two points at maximal intrinsic distance which are not vertices, the
result remains true, i.e. there are at least 5 geodesic segments connecting
the two points. Moreover, the convex polyhedra in the above class with two
points at maximal intrinsic distance connected by exactly 5 geodesic
segments form a dense set in the class. The results concerning the polyhedra
are not treated in this paper.

I would like to thank to Professor Tudor Zamfirescu for suggesting me the
problem for the case of convex surfaces, and for his constant support during
my stay at the University of Dortmund, Germany. My special thanks
to Professors  Jianguo Cao, Richard Hain, John Harer, David Morrison and
Fangyang Zheng for many useful comments and encouragement.

\vspace{0.2in}

{\bf 2. Preliminaries }

We recall some definitions. Complete explanations can be found in [2] and [4].

A function $f:I \to \R, I$ being an interval in $\R,$ is said to be
{\sl convex} if for $a\neq b, a, b \in I$ and $t \in (0,1)$ we have the
inequality
$$f(ta+(1-t)b)\leq tf(a)+(1-t)f(b)$$ The function is called {\sl strictly
convex} if the inequality is strict. If $N$ is a Riemannian manifold of
dimension $n \geq 2$, consider a geodesic segment parametrized proportionally
to arc length $\gamma :[0,1] \to N, \gamma (0)=p_0, \gamma (1)=p_1$. Denote by
$[p_0, p_1]$ the set $\gamma([0,1])$. A subset $V \subset N$ is called
{\sl convex}, if for every $p, q \in V$ there is a unique geodesic segment
from $p$ to $q$ and this is contained in $V$. A function $f:V \to \R \ $
is called {\sl (strictly) convex} if for every nontrivial geodesic
$\gamma : [0,1] \to V$ parametrized proportionally to arc length the function
$f\circ \gamma $ is (strictly) convex.

We now introduce a notion useful in what follows. For an open nonempty
convex set $V \subset N$ and for a point $p \in  V$, a convex function
$f: V  \to \R$ is called {\sl almost strictly convex} at $p$, if there
exists a line $l_p^f \subset \hbox{ T}_p V$ (passing through the origin)
such that
for every geodesic $\gamma : [0,1] \to V$ parametrized proportionally
to arc length  with $\gamma(0)=p$, $\dot \gamma(0) \not\in
l^f_p$, the function $f\circ \gamma$ is not constant. If $f$ is almost
strictly convex for all $p\in V$, we say that $f$ is almost strictly convex
on $V$. Of course, every strictly convex function is almost strictly
convex.

\vspace{0.1in}

{\bf Remark 2.1.} If $M_1$ and $M_2$ are Riemannian manifolds, then the
product Riemannian metric on $M_1 \times M_2$ is given by the action of
the two metrics on the product tangent space. So a curve

$$\gamma : [0,1] \to M_1 \times M_2, \quad \gamma(t)=(\gamma_1
(t),\gamma_2(t)),$$
is a geodesic parametrized proportionally to arc length in
$M_1\times M_2$ if and only if
$\gamma_i:[0,1]\to M_i, i=1,2$, are also geodesics parametrized proportionally
to arc length. It follows that if $U_i\subset M_i$ are nonempty convex sets
($i=1,2$), then $U_1\times U_2$ is also a convex set in $M_1\times M_2$.

\vspace{0.1in}

Consider now a simply connected complete Riemannian manifold $N$ of dimension
$n\geq 2$, with all the sectional curvatures bounded above by $\chi \leq 0$
(shortly, of non-positive curvature). We will say that $N$ has negative
curvature if $\chi < 0$.

We introduce some notations. For two distinct points $p_0, p_1 \in N$ consider
the geodesic segment $\gamma :[0,1] \to N$ parametrized proportionally to arc
 length, such that $\gamma (0)=p_0, \gamma (1)=p_1$, and for $t\ (0\leq t \leq
1)$ denote by $(1-t)p_0+tp_1$ the point $\gamma (t)$. Suppose that
$\Hi_\chi$ is the $2$-dimensional space of constant curvature
$\chi \leq 0$ (i.e. the euclidian plane for $\chi =0$, and the hyperbolic
plane with constant curvature $\chi$ for $\chi < 0$). For a triple
of points $[p_1, p_2, p_3]$ in $N$ the corresponding triple of points in
$\Hi_\chi$ is denoted by $[p^*_1, p^*_2, p^*_3]$, and has the property:
$$d(p_i, p_j)= d(p^*_i, p^*_j), 1\leq i<j\leq 3,$$
where $d^*$ is the distance function on $\Hi _\chi$. For a point $q\in
[p_1, p_2]$, the corresponding point $q^*\in [p^*_1, p^*_2]$ has the property
$d(p_1,q)= d^*(p^*_1, q^*)$.

It is known then that $N$ has non-positive (negative, respectively) curvature
in the metric sense, too
(see for example [4], page 52). This means that for every triple of points
$[p_1, p_2, p_3]$ in $N$ and the points $q\in [p_1, p_2], r\in [p_1, p_3]$
the corresponding triple $[p^*_1, p^*_2, p^*_3]$ and the points $q^* \in
[p^*_1, p^*_2], r^*\in [p^*_1, p^*_3]$ in $\Hi_\chi$ have the property:
$d(q,r)\leq d^*(q^*,r^*)$.

\vspace{0.1in}

{\bf Remark 2.2.} If $\Hi ^n $ is the $n$-dimensional hyperbolic space, then
the distance function $d:\Hi^n \times \Hi^n \to \R \ $ has the following
property:
$$d({p_1+q_1 \over 2}, {p_2+q_2 \over 2})< {d(p_1,p_2)+d(q_1,q_2) \over 2}$$
where $p_1, p_2, q_1, q_2$ are points in $\Hi ^n$ which are not on the same
maximal geodesic (see [3], page 37).

It is elementary that a continous  midconvex function is convex. It is also
easy to see that strict inequality for the middle of a segment implies that
the inequality is strict for every interior point of the segment. It follows
that, for $t \ (0<t<1)$, we have
$$d(tp_1+(1-t)q_1, tp_2+(1-t)q_2) \leq td(p_1, p_2)+ (1-t)d(q_1, q_2).$$

Of course, when the four points are on the same geodesic, the inequality
becomes equality. We need only the case $n=2$, and this can be also verified
by direct computation using the cosine formula in the hyperbolic plane. In
fact, the following is true.

\vspace{0.1in}

{\bf Lemma 2.3.} {\it Let N be a simply connected complete Riemannian
manifold of negative curvature as above, and $U_1$ and $U_2$ two nonempty
convex sets in $N$, such that $U_1 \cap U_2=\emptyset$. Then the restriction
of the distance function $d: U_1\times U_2 \to \R \ $ is almost strictly
convex.}

\vspace{0.1in}

P r o o f . From the remark 2.1., it follows that the set $U_1 \times U_2$
is a nonempty convex subset of the product manifold $N\times N$, so it makes
sense to talk about the convexity of the function $d: U_1\times U_2 \to
\R \ $. Take the points $p_1, q_1 \in U_1, p_2, q_2 \in U_2$, $p_i \not=
q_i, i=1,2$. Consider $t$, $0<t<1$; we shall prove that:
$$d(tp_1+(1-t)q_1, tp_2+(1-t)q_2)\leq td(p_1,p_2)+(1-t)d(q_1,q_2),$$
with equality only when $p_1, q_1, p_2, q_2$ are on the same geodesic.

Suppose they are not, and denote by $r_1$ the point $tp_1+(1-t)q_1$ and
by $r_2$ the point $tp_2+(1-t)q_2$. Consider also the point $r \in [p_1,
q_2]$, $r=tp_1+(1-t)q_2$. Using the comparison triangle $[p^*_1, q^*_1,
q^*_2]$ and the corresponding points $r^*_1, r^* \in \Hi_\chi$, we
obtain that:
$$d(r_1,r)\leq d^*(r^*_1, r^*) < (1-t) d^*(q^*_1, q^*_2)=(1-t) d(q_1, q_2),$$
where we have  assumed  that $q_2 \not\in [p_1, q_1]$ (possible, because
$p_1, q_1, p_2, q_2$ are not all on the same geodesic), and we have used
the remark 2.2. In a similar way, we have that:
$$ d(r, r_2)\leq t d(p_1, p_2),$$
and then, by addition that:
$$d(r_1, r_2)\leq d(r_1, r)+ d(r, r_2) < t d(p_1, p_2) + (1-t) d(q_1, q_2).$$
When the four points are on the same geodesic the inequality becomes
equality.  Consider now  a geodesic $\gamma :[-\delta, 1+\delta] \to N$,
 $\gamma (0)=p_1, \gamma (1)=p_2, \delta >0$;  from the inequality proved
above it can be inferred that the function $\quad
d: U_1\times U_2 \to \R \ $ is
convex. Furthermore, for a geodesic $ \quad \tilde \gamma : [0,1] \to N
\times N$ such
that $ \tilde\gamma(0)=(p_1, p_2)$, it follows that the function
$d\circ \tilde \gamma :[0,1] \to \R \ $ is strictly convex if and only if
$\displaystyle {\dot {\tilde \gamma}}(0)$, the tangent vector at $p$, is
not in the 2-plane
generated in T$_{(p_1,p_2)}(N\times N)$ by the vectors $(\dot \gamma(0), 0)$
and $(0, \dot \gamma(1))$. It's not difficult to see that the function
$d\circ \tilde \gamma :[0,1] \to \R \ $ is constant if and only if
the tangent vector $\displaystyle {\dot {\tilde \gamma}}(0)$ is
on the line given by the
vector $(\dot \gamma(0),\dot \gamma(1))$ in the tangent space to the
product manifold $N \times N$ at $(p_1, p_2)$. So the restriction
of the distance function to $U_1\times U_2$ is almost strictly convex and
the direction of the line $l^{\ d}_{(p_1, p_2)}$ (see the definition of an
almost strictly convex function) is given by the vector
$(\dot \gamma(0),\dot \gamma(1))$.

\vspace{0.1in}

{\it Note.} As Professor Jianguo Cao pointed out to me, lemma 2.3. is a
consequence of the proposition 1 from the paper of Schoen and Yau [6].
Our notion of almost strictly convexity for a function $f: V \to \R$
corresponds to the fact that the rank of the hessian of $f$ at the
considered point is at least $\hbox{dim}V -1$.

\vspace{0.1in}

For the case of non-positive curvature we will use the statement:

\vspace{0.1in}

{\bf Lemma 2.4.} {\it Let N be a simply connected complete Riemannian
manifold of non-positive curvature, $p \in N$ a fixed point, and $U_2$ a
 nonempty convex set in $N$ such that $p \notin U_2 $. Then the function
$f: U_2 \to \R \ $, defined by $f(q)=d(q,p)$ is almost strictly convex.
Moreover,
for any point $q\in U_2$, and any geodesic $\gamma : [0,1] \to N$,
 $\gamma(0)=q$, the function $f\circ \gamma $ is not constant (i.e. the
direction $l_q^f$ can be chosen arbitrarily).}

\vspace{0.1in}

P r o o f. In this case, the restriction of the distance function is
considered with one argument fixed. Using the same notations as in the
previous lemma, but keeping $p_1=q_1=p$, and by comparing
this time with triangles in the euclidian plane,
we have that:
$$d(p, r_2)\leq t d(p, p_2) + (1-t) d(p, q_2),$$
with equality if and only if the three points are on the same geodesic.
This shows that $f$ is almost strictly convex. Moreover, at a closer look,
for a point $q\in U_2$, even for the geodesic $\gamma :[0, 1+\delta]
 \to N$, $\gamma (0)=p, \gamma (1)=q, \delta >0$, the function
$f \circ \gamma$ (where the composition makes sense) is not constant. This
proves the last assertion in the lemma 2.4. Notice though that, in this
situation, the function $f$ falls short of being strictly convex -
the strict inequality does not hold in the direction $\dot \gamma (1)$, but
holds in all the others.

\vspace{0.1in}

We will need also the following elementary fact:

\vspace{0.1in}

{\bf Lemma 2.5.} {\it Consider in $\R^n (n \geq 2)$ for every i, $1\leq i
\leq k, k\leq n$, an ($n-1$)-dimensional linear subspace $H_i$, and denote
by $\overline S_i$ one of the closed halfspaces determined by $H_i$. If
$\displaystyle \bigcup^k_{i=1} \overline S_i = \R^n$, then
$\displaystyle \bigcap^k_{i=1} H_i$ is a linear subspace of dimension
at least $(n-k+1)$.}

\vspace{0.1in}

P r o o f. The lemma is clearly true if there exists $i,j$
($1\leq i \leq k$), such that $H_i=H_j$. Suppose that all the hyperplanes
are mutually distinct.

We use induction relative to $n$. For $n=2,3$ the lemma is true. Suppose
it's true for $n$. Then in $\R^{n+1}$, take $S'_i=S_i\cap H_1$, for every
$i$, $2\leq i \leq k$. Since $\displaystyle \bigcup^k_{i=1} \overline
 S_i = \R^n$,
and $H_1\not= H_i$, for all $i$, $2\leq i \leq k$, it follows that
$\displaystyle \bigcup^k_{i=2} \overline S'_i = H_1$. Using the induction we
have that:
$$\hbox{dim}(\bigcap^k_{i=2} H_i)\geq n-(k-1)+1=(n+1)-k+1,$$ so
$$\hbox{dim}(\bigcap^k_{i=1} H_i)\geq (n+1)-k+1,$$
which ends the proof of the lemma.

\vspace{0.1in}

{\bf 3. Main results}

The main tool in proving the theorems will be the following:

\vspace{0.1in}

{\bf Proposition 3.1.} {\it Suppose N is a Riemannian manifold of dimension
$m\geq 2$, not necessarily of negative curvature, $V$ an open convex set in
$N$, $p\in V$. Assume that there exists a natural number $k\geq 1$, such
that for every $i$, $1\leq i \leq k$, there is an almost strictly convex
function $f_i : V \to \R, f_i(p)=0$, with the property:

(i) there exists $\epsilon >0$, such that for every point $q\in B(p,\epsilon)
\subset V,$
$$\min _{1\leq i \leq k} f_i(q) \leq 0.$$
Then $k \geq m$.

Moreover, if the following condition holds also:

(ii) for every $i_1, i_2, \dots , i_m
(1 \leq i_1 < i_2 < \dots < i_m \leq k),$
$$ l^{f_{i_1}}_p \cap l^{f_{i_2}}_p \cap \dots \cap l^{f_{i_m}}_p = \{ 0 \},$$
then $k \geq m+1$.}

\vspace{0.1in}

P r o o f. We can suppose that $\overline {B(p,\epsilon)} \i V$ and
$B(p,\epsilon)$
 is convex (make $\epsilon$ smaller, if necessary). For every $i$, $1 \leq i
\leq k$, consider the open set
$$ A_i=\{ q \mid q \in B(p,\epsilon), f_i(q)< 0 \},$$ and the closed set
$$ B_i=\{ q \mid q \in \overline {B(p,\epsilon)}, f_i(q)\leq  0 \},$$

Of course, $\overline {A_i} \i B_i$, and since the functions $f_i$ are
convex, it follows that $A_i$ and $B_i$ are convex sets. For every $i$,
$ 1 \leq i \leq k$, consider the geodesic $\gamma _i : [-\epsilon, \epsilon]
\to \overline {B(p,\epsilon)}$, parametrized by arc length, such that
$\gamma _i(0)=p$ and $\dot \gamma _i(0) \in l^i_p$. Denote by $C_i$ the
set $\gamma_i ([-\epsilon,\epsilon])$. We make the following remark:

$$(3.2)\qquad \qquad \qquad \qquad  B_i \setminus \overline {A_i} \i C_i
\qquad \qquad$$

Indeed, if $B_i \setminus \overline {A_i} =\emptyset $, there is
nothing to prove. Suppose that $B_i \setminus \overline {A_i} \not=
\emptyset $. For every point $q \in B_i \setminus \overline {A_i}$, we have
that $f_i(q)=0$.

Consider first the case $A_i= \emptyset$. If $B_i$ consists of the single
point $p$, the property is clear. For another point $q \in B_i$, we have
that
$$f_i(q)=f_i(p)=0, [p,q] \i B_i, \hbox{and} f_i \geq 0 \ \hbox{on} \  \overline
{B(p,\epsilon)},$$
because $A_i=\emptyset$. The convexity of the function $f_i$ implies then
 that for every point $r\in [p,q], f_i(r)=0$, so $f_i$ is constant along
the segment $[p,q]$ and this gives $B_i \i C_i$.

Suppose next that there exists a point $p'\in A_i$. Then $f_i(p')<0$, and, in
fact, from the convexity of the function $f_i$ it follows that $f_i(p'')<0$,
for every $p''\in [p,p'], p'' \not= p$. So $p'' \in A_i$, which implies
that $p \in \overline {A_i}$. For a point $q \in B_i \setminus \overline
{A_i}$, we can find $\delta >0$ such that $B(q,\delta) \cap A_i = \emptyset$;
this means $f_i(r) \geq 0$, for every $r \in B(q,\delta)$. On the other
hand, the convexity of $f_i$ implies that for every point $r' \in [p,q]
\cap B(q,\delta), r' \not= q$, we have that $f_i(r') \leq 0$ (because
$f_i(p)=f_i(q)=0$). In conclusion $f_i(r')=0$, and $f_i$ is constant along
the segment $[p,q]$. So $[p,q] \i C_i$ , which ends the proof of the
relation 3.2.

Condition (i) in the proposition says in fact that
$$\displaystyle \bigcup^k_{i=1} B_i= \overline{B(p,\epsilon)}.$$
Combining this with the relation 3.2, we obtain that
$$\displaystyle \bigcup^k_{i=1} (\overline {A_i}\cup C_i)=
\overline{B(p,\epsilon)},$$
so
$$\displaystyle \overline{B(p,\epsilon)}\setminus \bigcup^k_{i=1}
\overline {A_i} \i \bigcup^k_{i=1} C_i.$$
But the difference of the two sets in the first part of the inclusion
is an open set in $\displaystyle \overline{B(p,\epsilon)}$ and
since $m= \hbox{dim}N \geq 2$, it is clear that the union of the geodesics
$C_i$ cannot cover this open set, unless
 the open set is empty,
$$\displaystyle \overline{B(p,\epsilon)}\setminus \bigcup^k_{i=1}
\overline {A_i}= \emptyset$$

Consider then $k_0$,
 such that $A_i\not= \emptyset$, for all $i$, $1\leq i
\leq k_0$, and $A_j=\emptyset$, for all $j$, $k_0 \leq j \leq k$. We will
prove that in fact $k_0 \geq m$ (respectively $k_0 \geq m+1$). Clearly we
have that
$$\displaystyle \bigcup^{k_0}_{i=1}
\overline {A_i}= \overline {B(p,\epsilon)}$$

For the non-empty convex set $A_i$, with $p \in \partial A_i$, we can apply
proposition 4.9.2. from [5]: there exists a support hyperplane $H_p^i \i
\hbox{T}_pN$, so that all the tangent vectors at $p$ to the geodesic
segments connecting  $p$ with points in the set $A_i$ are in the same
open half-space,
denoted by $S_p^i$, $1\leq i \leq k_0$, determined by $H_p^i$ in
$\hbox{T}_pN$.

The sets $\overline {A_i}$ can cover the closed ball $\overline {B(p,
\epsilon)}$ if and only if the corresponding closed half-spaces cover the
tangent space at $p$ (otherwise we would have a direction which is
in none of the closed half-spaces, so a geodesic segment which, at least
locally, is not contained in any of the  closed sets $\overline {A_i}$). So:
$$\displaystyle \bigcup^{k_0}_{i=1}\overline{S_p^i}=\hbox{T}_pN,$$
and then lemma 2.5. implies that:
$$\hbox{dim}\displaystyle (\bigcap^{k_0}_{i=1}H_p^i)\geq m-k_0+1.$$

Suppose now that $k_0 \leq m-1$; then:
$$\hbox{dim}\displaystyle (\bigcap^{k_0}_{i=1}H_p^i)\geq m-(m-1)+1=2$$
Define the set $U=B(p,\epsilon)\cap \hbox{exp}_p(\displaystyle \bigcap
^{k_0}_{i=1}H_p^i)$. Take a point $q \in U\setminus \{p\}$. Then there exists
$i_0, 1\leq i_0\leq k_0$, such that $q \in \overline {A_{i_0}}$. But from
the definition of the set $U$ we have that $q\not\in A_{i_0}$, so it
follows that $q\in \overline{A_{i_0}} \setminus A_{i_0}$ or in other
words $f_{i_0}(q)=0$. Using the convexity of the function and the fact that
$f_{i_0}(p)=0$, we get as above that $q\in C_{i_0}$. The  consequence of
this argument is that $\displaystyle{U\i \bigcup^k_{i=1}C_i}$. But this
is impossible, because $U$ is a submanifold of dimension at least 2, which
cannot be covered by finitely many  1-dimensional submanifolds. This ends
the proof of the first part of the proposition.

Consider now the case $k_0=m$ and suppose that the condition $(ii)$ holds.
Then
$$\hbox{dim}\displaystyle (\bigcap^{k_0}_{i=1}H_p^i)\geq 1.$$

Consider a vector $v\in \displaystyle{\bigcap^{k_0}_{i=1}H_p^i}$,
 $0<||v||<\epsilon$, and define the point $q=\hbox{exp}_p
v$. It is not possible that $q\in \displaystyle{\bigcup^{k_0}_{i=1}A_i}$,
 because this would mean that $v\in S_p^i$ for some $i$, which is not true.
On the other hand, if there is no $i$ such that $q\in \overline{A_i}$,
it would follow that there exists a neighborhood $W$ of $q$,
$W\i B(p,\epsilon)$,
 such that
$$W\cap \displaystyle (\bigcap^{k_0}_{i=1}\overline{A_i})=\emptyset,$$
which is impossible,too.

The argument shows that there exists $i_0$ such that
$q\in \overline{A_{i_0}}\setminus A_{i_0}$. The fact that $v$ is the
tangent vector at $p$ to the geodesic segment $[p,q]$ and the convexity
of the function
$f_{i_0}$ and of the set $\overline{A_{i_0}}$ imply that $[p,q]\i
\overline{A_{i_0}}\setminus A_{i_0}$ (no interior point of
 $[p,q]$ can be in $A_{i_0}$, and $f_{i_0}$ is constant
along this geodesic segment). This means that $[p,q] \i C_{i_0}$, so
in fact $v\in l_p^{f_{i_0}}$. If $k_1$ is the number of $A_i$'s with
the property that $q\in \overline{A_i}$, condition $(ii)$ gives that
$k_1\leq m-1$. The first part of the proposition can be applied now, and a
contradiction is obtained. This ends the proof of the proposition.

\vspace{0.1in}

{\bf Remark 3.3.} Notice that, when conditions $(i)$ and $(ii)$ are satisfied,
the argument from the last part of the proof shows that there exists
a neighborhood $W$ of $p$ such that for $q\in W, q\not= p$, we have
$$\min _{1\leq i \leq k} f_i(q) < 0$$
In other words, the function $\displaystyle \min _{1\leq i \leq k} f_i$
has a strict local maximum at $p$.

\vspace{0.1in}

We are now in the position to prove the main results.

\vspace{0.1in}

{\bf Theorem 3.4.} {\it Let M be a complete Riemannian manifold of
dimension $n\geq 2$ and negative curvature (all the sectional
curvatures bounded above by $\chi <0$). If $(p_1, p_2)\in M\times M$ is
a local maximum for the distance function on $M$, $d:M\times M\to \R$, then
the points $p_1$ and $p_2$ are connected by at least $2n+1$ distinct geodesic
segments.}

\vspace{0.1in}

P r o o f. Consider the universal covering space $\wm$ which is
diffeomorphic to $\R^n$ and $\pi : \wm \to M$ the covering map. Denote
by $\tilde d$ the distance function on $\wm$, and by $d_0$ the distance
between the points $p_1$ and $p_2$. Take a point $\tilde p_1 \in
\wm$ so
that $\pi (\tilde p_1)=p_1$.

The closed ball $\displaystyle{\overline{B(\tilde p_1,d_0)}}$ is compact in
$\wm$. From the discreteness of the fibre $\pi^{-1}(p_2)$, it follows
that $p_1$ and $p_2$ are connected in $M$ by just finitely many geoedsic
segments
$$\gamma_i:[0,1] \to M, 1\leq i \leq k, \gamma_i(0)=p_1,
\gamma_i(1)=p_2,$$
parametrized proportionally to arc length. Consider their
liftings to the universal cover
$$\tilde \gamma_i:[0,1] \to \wm, 1\leq i \leq k, \tilde \gamma_i(0)=
\tilde p_1,  \tilde \gamma_i(1)=\tilde p^i_2 \in \pi^{-1}(p_2),$$
The discreteness of the fibre implies also that it is possible to find
$\epsilon _0>0$, so that $\tilde d(\tilde p_1, \tilde p)> d_0 + \epsilon_0$,
for every $\tilde p \in \pi ^{-1}(p_2) \setminus \{ \tp ^1_2,\dots ,
\tp^k_2 \}$. Define the sets:
$$\displaystyle U_1=B(p_1,{\epsilon_0 \over 4}), U_2=B(p_2,
{\epsilon_0 \over 4}),$$
$$\displaystyle \wu_1=B(\tp_1,{\epsilon_0 \over 4}), \wu_2^i =B(\tp_2^i ,
{\epsilon_0 \over 4}), 1\leq i\leq k.$$

We can suppose that $\epsilon_0$ is small enough so that $U_1 \cap U_2=
\emptyset$, $\wu_2^i \cap \wu_2^j =\emptyset$, for $i\not= j$, and the
restrictions $\pi_1:\wu_1 \to U_1$, and $\pi_{2,i}: \wu_2^i \to U_2$ of
the covering map are isometries. The negative curvature and the convexity
of the distance function imply that the balls $U_1, U_2, \wu_1, \wu_2^i$ are
convex.

For every $i$, $1\leq i \leq k$, consider the restriction of the distance
function $\tilde d$ on $\wm$, $\tilde d_i:\wu _1 \times \wu_2^i \to \R$, and
the map $\pi_i : \wu _1 \times \wu_2^i \to U_1 \times U_2$, defined by
$\pi_i=(\pi_1, \pi_{2,i})$. We have that $\pi_1$ and $\pi_{2,i}$ are
isometries, so $\pi_i$ is also an isometry. Define the function
$f_i: U_1 \times U_2 \to \R$, by the relation $f_i=\tilde d_i \circ \pi_i
^{-1}$. In fact, for $q_1\in U_1, q_2 \in U_2$, it's clear that
$f_i(q_1,q_2)=\tilde d (\tilde q_1, \tilde q_2^i)$, where
$\tilde q_1=\pi _1^{-1}(q_1),$ and $\tilde q_2^i=\pi _{2,i}^{-1}(q_2).$
The claim is that:
$$(3.5)\qquad \qquad \qquad  \displaystyle d(q_1,q_2)=
\min_{1\leq i \leq k} f_i(q_1,q_2), \quad \hbox{for} \quad q_1\in U_1,
q_2 \in U_2 \qquad$$

For, notice that $\displaystyle d(q_1,q_2)=\min_{\tilde q_2 \in \pi ^{-1}
(q_2)} \tilde d (\tilde q_1,\tilde q_2)$. But:
$$\displaystyle \tilde d (\tilde q_1,\tilde q_2^i) \leq \tilde d (\tilde q_1,
\tilde p_1)+ \tilde d(\tp_1,\tp_2^i)+\tilde d(\tp _2^i,\tilde q_2^i) <
{\epsilon \over 4} + d_0 + {\epsilon \over 4} = d_0 + {\epsilon \over 2}$$
If $\tilde q_2 \in \pi ^{-1}(q_2) \setminus \{ \tilde q ^1_2,\dots ,
\tilde q ^k_2 \}$, then:
$$\displaystyle \tilde d (\tilde q_1,\tilde q_2) \geq \tilde d(\tp_1,\tp)-
\tilde d(\tp, \tilde q_2)- \tilde d(\tilde q_1, \tp_1) > d_0 +\epsilon
-{\epsilon \over 4}-{\epsilon \over 4}=d_0+{\epsilon \over 2},$$
where $\tp$ is the point of the fibre $\pi^{-1}(p_2)$ which is closest
to $\tilde q_2$. The relation $(3.5)$ is proved.

We have that $\wu_2^i \cap \wu_1=\emptyset$, because $U_1 \cap U_2=\emptyset$.
By applying the lemma 2.3, we obtain that the functions $\tilde d_i$ are
almost strictly convex, and that the lines $l_{\tp_i}^{\tilde d_i}$,
$\tp_i=(\tp_1,\tp_2^i)$, are generated by the vectors
$(\dot {\tilde \gamma}_i(0), \dot {\tilde \gamma}_i(1))$. From the
equality
$$f_i \circ \gamma_i=(f_i \circ \pi_i)\circ (\pi_i^{-1}\circ \gamma_i)=
\tilde d_i \circ \tilde \gamma_i,$$
and using the fact that $\pi _i$ are isometries, it follows that the
functions $f_i$ are almost strictly convex with $l_p^{f_i}= (\pi_i)_*
(l_{\tp_i}^{\tilde d_i})$, $p=(p_1,p_2)$. In order to apply the proposition
3.1 for the functions $f_i$, we have to prove that the conditions $(i)$
and $(ii)$ are satisfied.

The fact that $p=(p_1,p_2)\in M \times M$ is a local maximum for the
restriction to $U_1 \times U_2$ of the distance function on $d: M\times M
\to \R$, together with the equality 3.5, provide exactly the
condition $(i)$ for
the functions $f_i$. This implies that $k\geq 2n\geq 4$. Next, we prove
that, in this particular case, condition $(ii)$ holds in the stronger
form:

$(ii')\quad \hbox{for every} \quad i_1, i_2, i_3 \
(1\leq i_1 \leq i_2 \leq i_3 \leq k),$
$$l_p^{f_{i_1}}\cap l_p^{f_{i_2}}\cap l_p^{f_{i_3}}=\{ 0 \}$$

For, suppose that there exists a non-zero vector $v=(v_1,v_2)$,
$v\in l_p^{f_{i_1}}
\cap l_p^{f_{i_2}}\cap l_p^{f_{i_3}}$. But the direction of the
line $l_p^{f_i}$, $i \in \{ i_1, i_2, i_3 \}$, is given by the vector
$(\dot \gamma_i(0), \dot \gamma_i(1))$. Since $v\not= 0$, one of its
components is non-zero, too. Suppose $v_1\not= 0$. It follows that
$v_1=\pm \dot \gamma_i(0)$, for all  $i \in \{ i_1, i_2, i_3 \}$ (the
geodesic segments joining $p_1$ and $p_2$ are parametrized proportionally
to arc length and they have the same length). But this
implies that at least two of the  considered geodesic segments have
the same tangent vectors at $p_1$, which is impossible. The contradiction
obtained shows that the condition $(ii')$ holds for the functions $f_i$. The
proposition 3.1 can be applied in this case, which ends the proof of the
theorem.

\vspace{0.1in}

{\bf Remark 3.6.} If we translate the remark 3.3 in the language of the
theorem 3.4, it follows that the  pairs of points, for
which the distance function has a local maximum, are isolated in
the product topology.

\vspace{0.1in}

{\bf Theorem 3.7.} {\it Let M be a complete Riemannian manifold of
dimension $n\geq 2$ and non-positive curvature (all the sectional
curvatures bounded above by $\chi \leq 0$), $p_1$ a fixed point in $M$. If
$p_2 \in M$ is a local maximum for the function $f: M \to \R$, defined
 by $f(q)=d(q,p_1)$, ($d$ is the distance function on $M$), then
the points $p_1$ and $p_2$ are connected by at least $n+1$
distinct geodesic segments.}

\vspace{0.1in}

P r o o f. We will use the notations from the proof of the
previous theorem. The universal covering space $\wm$ is again
diffeomorphic to $\R^n$, and using the same construction we
define the functions $f_i :  U_2 \to \R$ by the relation
$f_i (q_2)=\tilde d (\tilde p_1, \tilde q_2^i)$. The same argument
as above gives that $\displaystyle f(q_2)= d(p_1,q_2)=
\min_{1\leq i \leq k} f_i(q_2), \quad \hbox{for} \quad
q_2 \in U_2 $. Notice that in this case, lemma 2.4. implies that the
functions $f_i$ are almost strictly convex, but also that the ``bad"
directions $l_{p_2}^{f_i}$ can be chosen arbitrarily. This means that
condition $(ii)$ from the proposition 3.1. is satisfied easily with
a general choice of the directions $l_{p_2}^{f_i}$. Proposition 3.1.
implies $k \geq n+1$, which ends the proof of the theorem.

\vspace{0.1in}

{\bf Corollary 3.8.} {\it On every compact Riemannian manifold of negative
curvature, there are pairs of points connected by at least $2n+1$
distinct geodesic segments (for example, the points at maximal
distance).}

\vspace{0.1in}

{\bf Corollary 3.9.} {\it On every compact Riemannian manifold of non-positive
curvature, for any given point $p_1 \in M$,  there exist a point
$p_2 \in M$ (for example the point at maximal distance from $p_1$) and at
least $n+1$ distinct geodesic segments connecting the two points.}

\vspace{0.1in}

{\bf Remark 3.10.} The flat torus example given in the introduction shows
that corollary 3.9. is the best one can expect in general for the case
of non-positive curvature. But the flatness is a very restrictive situation.
It would be interesting to find conditions in which the result can be
improved (for example, if we set as a hypothesis that the two points
which realize the diameter of the manifold are isolated in the product
topology).

\vspace{0.1in}

{\bf Remark 3.11.} As one can see, the definitions and the proofs, except
the proof of proposition 3.1. are essentially metric. It is likely that
the results are true without the differentiability hypothesis for
spaces of non-positive curvature (see [1], [4]) which
are $n$-dimensional topological manifolds. What one will need is the
generalization of the tangent space for metric spaces, the so-called
space of directions. The notion is discussed in [1], but the details of
the proof of the analog of the proposition 3.1 seem to be more difficult. A
case which probably can be studied directly is the case of polyhedra of
non-positive curvature, where the same methods as above should work.

\vspace{0.2in}

{\bf References}

\vspace{0.1in}

[1] A.D. Alexandrov, V.N. Berestrovskii \& I.G. Nikolaev, {\it Generalized
Riemannian spaces}, Russian Math. Surveys {\bf 41} (1986) 1-54.

[2] W. Ballmann, M. Gromov \& V.Schroeder, {\it Manifolds of nonpositive
curvature}, Progress in Math., Vol. 61, Birkh\"auser, Boston, 1985.

[3] R. Benedetti \& C. Petronio, {\it Lectures in hyperbolic manifolds},
Universitext, Springer, Berlin, 1992.

[4] E. Ghys \& P. de la Harpe (Ed.), {\it Sur les groupes hyperboliques
d'apr\`es Mikhael Gromov}, Progress in Math., Vol. 83, Birkh\"auser, Boston,
1990.

[5] H. Karcher, {\it Schnittort und konvexe Mengen in vollst\"andigen
Riemannschen Mannigfaltigkeiten}, Math. Ann. {\bf 177} (1968) 105-121.

[6] R. Schoen \& S.T. Yau, {\it Compact group actions and the topology of
manifolds with non-positive curvature}, Topology {\bf 18} (1979) 361-380.

\vspace{0.25in}

{\it Address:} Department of Mathematics, Duke University,
 Box 90320,
\newline{Durham, NC 27708-0320, U.S.A.}

{\it Current address (until May 1, 1995):} Department of
 Mathematics,
\newline{Cornell University, Ithaca, NY 14853, U.S.A.}

{\it E-mail:} horja$@@$math.duke.edu

\end{document}